\documentclass[11pt,reqno]{amsart}
\usepackage{amssymb,verbatim}
\usepackage[mathcal]{euscript}
\def\CC {{\mathbb C{\,}}} \def\RR {{\mathbb R{\,}}}

  \let\phi\varphi

\newcommand{\RE}{\operatorname{Re}} 
\DeclareMathOperator{\sign}{sgn}

\DeclareMathOperator{\Cr}{Cr}
\def\DOD#1#2{ \frac{\partial {#1}}{\partial {#2}} }

  \theoremstyle{definition}
  \newtheorem*{EXAM*}{Example}
  \newtheorem*{DEF*}{Definition}

  \theoremstyle{remark}
  \newtheorem*{NOTE*}{Note}
  \newtheorem*{REM*}{Remark}

  \theoremstyle{plain}
  \newtheorem*{COR*}{Corollary}
  \newtheorem*{LEM*}{Lemma}
  \newtheorem*{THM*}{Theorem}
  \newtheorem*{PROP*}{Proposition}

\theoremstyle{definition}
\newtheorem{DEF}{Definition}

\theoremstyle{remark}

\theoremstyle{plain}
\newtheorem{THM}{Theorem}
\newtheorem{PROP}{Proposition}

\begin{document}


\address{Graev, Mark Iosifovich. NIISI RAN. Nahimovskii prospect 36, коrp.~1. Moscow 117218  Russia}
\email{graev\_36@mtu-net.ru}
\address{Litvinov, Grigory Lazarevich. Nagornaya str., 27, korp.~4, kv.~72. 117186  Moscow  Russia}
\email{glitvinov@gmail.com}
\author{M.I.~Graev, G.L.~Litvinov}
\thanks{
The first author is supported by the RFBR grant 10-01-00041a.
}
\thanks{
The second author is supported by the joint RFBR (Russia) and CNRS (France) grant.
}
\title{Integral geometry, hypergroups,\\ and I.M. Gelfand's question }



\begin{abstract}This note is an attempt to give an answer for the following old I.M.~Gelfand's question: why some important problems of integral geometry (e.g., the Radon transform and others) are related to harmonic analysis on groups but for other quite similar problems such relations are not clear? In the note we examine standard problems of integral geometry generating harmonic analysis (the Plancherel theorem etc.) on pairs of commutative hypergroups that are in a duality of Pontryagin's type. As a result new meaningful examples of hypergroups are constructed.
 
\end{abstract}

\maketitle

\large 

The subject of integral geometry (in the sense of \cite{G,GGV,GGG})
is formed by integral transforms mapping functions on a manifold (or space)
$X$ to their integrals over submanifolds in $X$ forming a family  $\hat X$ of submanifolds in $X$, so new functions defined on $\hat X$ appear. One of basic problems in integral geometry is a reconstruction of the initial function on $X$ starting from its image on $\hat X$.

This note is an attempt to give an answer for the following old I.~M.~Gelfand's question: why some important problems of integral geometry (e.g., the Radon transform and others) are related to harmonic analysis on groups but for other quite similar problems such relations are not clear? In the note we examine standard problems of integral geometry generating harmonic analysis (the Plancherel theorem etc.) on pairs of commutative hypergroups that are in a duality of Pontryagin's. As a result new meaningful examples of hypergroups are constructed.


{\bf 1. Hypergroups and generalized Fourier transforms}.
Suppose that $X$ and $\hat X$  have structures of smooth manifolds or diffeological spaces with fixed measures $dx$ and $dy$ respectively.

Let us examine transforms of the form
\begin{equation}\label{1}
F: f(x)\mapsto \hat f(y)=\int f(x) e(x,y)\; dx,
\end{equation}
where $e(x,y)$ is a generalized function on $X\times\hat X$. Suppose that this transform generates an isomorphism  $L^2(X,dx)\to L^2(\hat X, dy)$, such that the following generalized {\it Plancherel formula} is valid:
\begin{equation}\label{2}
\int f(x){\overline{g(x)}} dx = \int \hat f(y){\overline
{\hat g(y)}} dy.
\end{equation}

In addition suppose that this isomorphism can be extended to
$\delta$-functions and set $g(x)=\delta_x(\cdot)$; then from
\eqref{2} it follows that $\hat{\delta}_x=e(x,y)$ and the following inversion formula is valid:
\begin{equation}\label{3}
f(x)=\int\hat f(y){\overline{e(x,y)}} dy.
\end{equation}

On the other hand the Plancherel formula  \eqref{2} follows from \eqref{1}
and \eqref{3}.
We suppose that in the space  $L^2 (X,dx)$ there is a dense locally convex linear subspace  $S$ (or lineal for the sake of brevity), such that  $S$ consists of continuous functions (but not all of them) and the lineal $\widehat S=F(S)$ possesses the same properties with respect to the space $L^2(\widehat X,dy).$ It is assumed that $S$ and
$\widehat S$ are algebras with respect to the usual multiplication of functions and the isomorphism \eqref{1} can be extended to  $\delta $-functions belonging to the spaces  $M$ and  $\widehat M$, dual to $S$ and $\widehat S$ respectively.

The described construction is a formalization and generalization of heuristic ideas of B.~M.~Levitan \cite{Le,LL}.

In the case at hand, when formulas \eqref{1} and \eqref{3} and hence \eqref{2} are valid, we shall say that the transform $F$ is a {\it  generalized Fourier transform}
(or GFT for the sake of brevity); functions $e(x,y)$ and $\overline{e(x,y)}$ will be called  {\it generalized exponential functions.}

\begin{PROP}{}
Under the specified assumptions,  $X$ has a structure of {\it commutative hypergroup} in the following sense: generalized translation operators act in  $S$ and for these operators the associativity axiom of J. Delsarte \cite{De} is valid; generally speaking, it is not supposed that X has a neutral element. These generalized translation operators are defined by the formula
\begin{equation}{}
R^y f(x)=\int \hat f(\chi) \cdot \overline{e(y,\chi )} \cdot\overline{e(x,\chi )} d\chi ,
\end{equation}
where $ \chi\in\widehat X$ defines the character $\chi(x)= e(x, \chi )$ on $X.$
Similarly
$x(\chi )=\overline{e(x, \chi )}$ is a character on $\widehat X.$
\end{PROP}

The hypergroup commutativity means that every two generalized translation operators commute. Then $\widehat X$ is a hypergroup dual to $X$;
of course, the hypergroup $X$ is dual to $\widehat X$ and
$F^{-1}(\delta _y)=\overline{e(x,y)}.$

For example, if $X=\widehat X=\RR^n,$ $dx$ and $dy$ are normalized invariant measures, then $F$ is the usual Fourier transform,
$ e (x,y)$ is the exponential function $e^{i{\langle x,y \rangle}},$
and for $S$ it is possible to use the L.~Schwartz space of functions
which are rapidly decreasing with all the derivatives. In this case the space $M$ consists of all tempered generalized functions (distributions).
Below in Section 5 the space $S$ is the space of all smooth functions on a
compact hypergroup $X$; in other sections $S$ is the L.~Schwartz space on a hypergroup $X$. 

Usually generalized translation operators act also in the spaces of smooth functions, summable functions, measures and generalized functions with compact supports etc. 
%
%
Different versions of the concept of hypergroup are discussed, e.g., in
\cite{De}--\cite{Li}. Of course, every group (or semigroup) $G$ is an example of a hypergroup; in this case translations $R^x$ are of the form:
$$
R^x: f(t) \mapsto f(tx),
$$
where $tx$ is the product of elements of the group.

There are quite few examples of nontrivial meaningful hypergroups related to harmonic analysis. The most well known one is the Gelfand
pair related to harmonic analysis of spherical functions, see, e.g., \cite{Li,LL}. One of our aims in this note is to extend the collection
of examples of this type.



{\bf 2. GFT associated with the generalized Radon transform on
$\CC^n$}.
The Radon transform on $\CC^n$ maps functions $f$ on
$\CC^n$ to their integrals over hyperplanes in $\CC^n,$ i.e. functions $\mathcal Rf$ on the manifold of hyperplanes in  $\CC^n.$ 
For an explicit description of this transform we shall define hyperplanes by the equations $x_n=\sum_{i=1}^{n-1}a_ix_i+a_n$
and we shall use $a=(a_1, \dots ,a_n)$ as local coordinates for the manifold of hyperplanes in  $\CC^n$.  Using these notations and the delta-functions notation we can render the Radon transform in the following form: 
$$
(\mathcal R f)(a)=\int_{\CC^n} f(x) \delta (x_n-a_1x_1- \ldots
-a_{n-1}x_{n-1}-a_n) d\mu(x),
$$
where $\delta (\cdot )$ is the delta-function on $\CC,$ and $d\mu(x)$ is the Lebesgue measure on $\CC^n.$ It is known that the following inversion formula is valid:
if $\phi=\mathcal R f,$ then
$$
f(x)=\int_{\CC^n}\phi(a) \delta^{(n-1,n-1)} (x_n-a_1x_1- \ldots
a_{n-1}x_{n-1}-a_n) d\mu(a),
$$
where
$
\delta^{(n-1,n-1)}(t)
= \DOD{^{2n-2}}{t^{n-1} \partial \overline t^{n-1} } \delta (t),
$

We shall say that a {\it generalized Radon transform} on $\CC^n$
associated with an arbitrary generalized function $u(t)$ on $\CC$ is
the integral transform
$$
(\mathcal R_u f)(a)=\int_{L^n} f(x) u (x_n-a_1x_1- \ldots -a_{n-1}x_{n-1}-a_n)
d\mu(x).
$$
\begin{THM}{}\label{THM:1} If the Fourier transform $\widetilde u(c)$ of 
$u(t)$ is an usual function and this function is nonzero almost everywhere, then the generalized Radon transform $\mathcal R_u$ is invertible and the inversion formula is of the form
\footnote{Here and in the next integral formulas we suggest that coefficients which usually are placed before the integral symbol are
included in the corresponding measure normalization.
}:
if $\phi=J_u f,$ then
\begin{equation}{}\label{3-71}
f(x)=\int_{\CC^n} \phi(a) U(x_n-a_1x_1- \ldots   a_{n-1}x_{n-1}-a_n) d\mu(a),
\end{equation}
where  $U(t)=\int_{\CC} [\widetilde u(c)]^{-1} |c|^{2(n-1)} e^{i\RE (ct)}d\mu(c).$
\end{THM}

\begin{COR*}{}
If $|\widetilde u(c)|=|c |^{n-1},$ then the generalized Radon transform $J_u$ is a GFT, i.е. $U(t)=\overline{u(t)}.$
\end{COR*}

{\it Examples of GFT. } Functions $\widetilde u(c)=c^k\overline c^{n-k-1},$ $k=0,1, \dots ,n-1$ correspond to
local GFT with kernels $u(t)= \delta ^{(k,n-k-1)}(t).$
Functions $\widetilde u(c)=c^{\lambda }\overline c^{\mu},$ which differ from functions indicated above, where $\RE(\lambda +\mu)=n-1$ and  $\lambda -\mu$ is an integer number, correspond  to GFT with nonlocal kernels $u(t)=t^{- \lambda -1}\overline t^{-\mu-1}.$

{\bf 3. GFT associated with generalized Radon transform on $\RR^n$}.
Definitions of the usual (see \cite{Ra}) and generalized Radon transform on $\RR^n$ are similar to the complex case. In the real case the generalized Radon transform  $\mathcal R$ is a GFT if the Fourier transform  $\widetilde u(c)$ of the kernel $u(t)$ satisfies the equation $|\widetilde u(c)|=|c|^{\frac{n-1}{2}}.$
Hence (see \cite{G-S})we have one and only one local GFT for odd $n,$ $n=2k+1$ only. This kernel corresponds to the function 
$\widetilde u(c)=c^k,$ and its kernel has the form $\delta ^k(t).$  Examples of GFT with nonlocal kernels are GFT with kernels $u(t)=|t|^{-\frac{n+1}{2}-i\rho},$
$\rho\in\RR.$

\begin{REM*}{} Similarly it is possible to construct GFT for  generalized Radon transforms on $L^n,$
where $L$ is a continuos non-Archimedean locally compact field.
\end{REM*}

{\bf 4. GFT associated with a transform of functions on the sphere ${S^n \subset \RR^{n+1}}$}.
Let us examine an integral transform mapping {\it even} functions on the sphere $S^n$ to their integrals on geodesic hypersurfaces (analogs of big circles on  $S^2$). In the spherical coordinates 
$\omega $ these hypersurfaces are defined by the equations ${\langle \xi, \omega  \rangle}{=}0,$ and the integral transform can be presented in the form
\begin{equation}{}\label{8-5}
f(\omega ) \to (Jf)(\xi)=\int_{S^n} f(\omega )
\delta ({\langle \xi, \omega  \rangle}) d \omega ,
\end{equation}
where $d \omega $ is an invariant measure on the sphere. There exists an inversion formula representing the initial function $f$ via its image under the transform $Jf.$

For every $\lambda \in\CC$ let us define the following generalized transform for even functions on $S^n$ by the equation
\begin{equation}{}\label{18-5}
( J_{\lambda }f) (\xi)=\int_{S^n} f(\omega )
\frac{|{\langle \xi, \omega  \rangle}|^\lambda}
{\Gamma (\frac{\lambda +1}{2})}  d\omega .
\end{equation}
The integral converges for $\RE  \lambda >-1$ and it is defined for every  $\lambda \in\CC$ as an analytic continuation in $\lambda$. In particular, $J_{-1}=J.$

\begin{THM}{}\label{THM:3-7} The following inversion formula holds:
\begin{equation}{}\label{18-6}
\text{если} \quad  \phi=J_\lambda f,\quad\text{то}\quad f=J_{- \lambda -n-1}\phi .
\end{equation}
\end{THM}
\begin{COR*}{} The transform $J_\lambda $ is a GFT for $\lambda =-\frac{n+1}{2}+i\rho,$ $\rho\in\RR$.
\end{COR*}
\noindent
In particular, for $n=4k+1$ and $\lambda=-2k-1$ we have the GFT with the local kernel  $ \delta ^{2k}(\cdot ). $

Similarly define the transform $J_u$ for {\it odd}
functions on the sphere replacing $|{\langle \xi, \omega  \rangle}|^\lambda$  with $|{\langle \xi, \omega  \rangle}|^\lambda
\,\sign({\langle \xi, \omega  \rangle})$
 in \eqref{18-5}.
For this definition, the same inversion formula \eqref{18-6} holds and
the transform $J_u$ is a GFT under the same conditions as in the case of 
even functions. The difference is that local GFT exist in the case of $n=4k-1,$ and in this case the kernel is the odd generalized function $\delta ^{2k-1} (\cdot ).$  Note that our hypergroup is compact and the dual hypergroupm is compact too.

{\bf 5. GFT related to the complex of $k$-dimensional planes in  $\CC^n$}.
Define a family of $n$-dimensional submanifolds (complexes) in the manifold of all $k$-dimensional planes in the space $\CC^n,$ $0<k<n-1.$

Let the space $\CC^n$ be represented in the form of direct sum $\CC^n=\CC^k \oplus \CC^l,$ ${k+l=n}$ of the spaces with coordinates $x=(x_1, \dots, x_k)$ and
$y=(y_1, \dots ,y_l)$.  Define $k$-dimensional planes in $\CC^n$
by the following equations solved with respect to the coordinates $y_i:$
$y_i=u_{i1}x_1+ \ldots + u_{ik}x_k+ \alpha _i,$ $i=1, \dots ,l.$

Let us fix an arbitrary $(l,k)$-matrix $u(t)=\|u_{ij}(t)\|$ whose elements are polynomials of $t=(t_1, \dots ,t_k)\in\CC^k.$ 
For this matrix let us construct the submanifold of $k$-dimensional planes defined by the equations
\begin{equation}{}\label{6-421}
y_i=u_{i1}(t)x_1+ \ldots + u_{ik}(t)x_k+ \alpha _i ,
\quad  i=1, \dots ,l,
\end{equation}
or $y=u(t)x + \alpha$ for short.

The condition that the submanifold $K$ just defined is a complex, i.e.
$\dim K=n,$ is equivalent to the condition of nonsingularity for the map $t \mapsto u(t)$ of the space $\CC^k$ to the space of $(l,k)$-matrices.
We shall use vectors $\alpha =( \alpha _1, \dots ,\alpha _l )$ and
$t=(t_1, \dots ,t_k)$ as coordinates for $K.$

Complexes $K$ have a simple geometric structure. Specifically, equations  \eqref{6-421} for $\alpha =0$ define a $k$-dimensional family of $(k-1) $-dimensional planes on a hyperpane in $\CC^n$ at infinity. The complex $K$ consists of all $k$-dimensional planes containing at least one of these $(k-1)$-dimensional planes. In particular, for $k=1$ the complex $K$ consists of straight lines in  $\CC^n$ that intersec a fixed curve in a hyperplane at infinity.

The complex $K$ generates an integral transform $J$ mapping functions 
$f(x,y)$ on $\CC^n$ to their integrals over planes of the complex. Using the delta-functions notation, it is convenient to present this transform in the form
\begin{equation}{}\label{6-46}
(Jf)(\alpha ,t)=\int_{\CC^n}f(x,y) \,\delta (y- u(t)x- \alpha  ) \,d\mu(x) \,d\mu(y),
\end{equation}
where $\delta (\cdot )$ is the delta-function on $\CC^k,$  and $d\mu(x)$,  $d\mu(y)$ are Lebesgue measures on $\CC^k$ and $\CC^l.$

We now examine general transforms $J_a$ which are defined by replacing
in  definition \eqref{6-46}
the delta-functions $\delta (s)$ with arbitrary generalized functions $a(s,t)$ and find inversion formulas for these transforms.
By definition we have
\begin{equation}{}\label{6-52}
(J_af)(\alpha ,t)=\int f(x,y)\, a(y-u(t)x- \alpha ,t)\, d\mu(x)\,d\mu(y).
\end{equation}
Note that  $J_af$ is a result of convolution with respect to $\alpha $ of the function $\phi=Jf$ with the function $a(s,t)$ :
 \begin{equation}{}\label{6-53}
(J_af)( \alpha  ,t ) =\int_{\CC^l} \phi(\alpha +s,t)\,a(s,t)\, d\mu(s).
\end{equation}

There is a simple relation between the Fourier transform  $\widetilde f(\eta,\xi)$ of the function $f(x,y)$ and the Fourier transform 
$\widetilde\phi(\xi,t)$ with respect to $\alpha $ of the function $\phi=J_af:$
\begin{equation}{}\label{6-59}
\widetilde\phi(\xi,t)=\widetilde f(-\xi u(t),\,\xi) \widetilde a(-\xi,t),
\end{equation}
where $\widetilde a(\xi,t)$ is the Fourier transform with respect to $s $ of the function $a(s ,t).$

Let us check that {\it if $\widetilde a(\xi,t)$ is an usual function which is nonzero almost everywhere, then the function $f$ can be reconstructed from the function $\phi=J_af$} in a unique way.

Note that  by virtue of \eqref{6-59} this condition implies that the function 
$F(\eta,\xi)=\widetilde f(-\xi u(t),\xi)$ can be reconstructed from the function $\phi$ uniquely. On the other hand from the nonsingularity of $K$ and the analyticity condition it follows that for almost every pair $(\eta,\xi) $ there exists  $t\in\CC^k$ such that $-\xi u(t)=\eta.$
It then follows that the function  $\widetilde f,$ and hence the function $f$ can be reconstructed from the function $F.$ Let us describe an inversion formula explicitly.

\begin{DEF}{}\label{DEF:5-1} We shall say that {\it the Crofton function}
related to the complex $K,$ is the function $\Cr_K(\eta,\xi)$ on $\CC^k \oplus \CC^l,$  equal to the number of solutions  $t$ of the equation $\eta=-\xi u(t).$ If the Crofton function is constant almost everywhere, then this constant is called {\it the Crofton number}. Denote this number by $\Cr_K.$
\end{DEF}
From the condition on the complex $K$ it follows that its Crofton function is nonzero and constant almost everywhere.

\begin{THM}{}\label{THM:6-9} If $\widetilde a(\xi,t)$  is a usual function which is nonzero almost everywhere, then there is the following inversion formula for the integral transform $f \mapsto \phi=J_af:$
\begin{equation}{}\label{6-55}
f(x,y)=\int \phi(\alpha ,t )\,A(y-u(t)x- \alpha ,t)\,d\mu(\alpha )\,d\mu(t),
\end{equation}
where
\begin{equation}{}\label{6-56}
A(s,t)=\frac{1}{\Cr_K} \int [\widetilde a(\xi,t)]^{-1}
\Bigl|\DOD{(\xi u(t))}{t}\Bigr|^2
 e^{i \RE  {\langle s,\xi \rangle}} d\mu(\xi).
\end{equation}
\end{THM}

\begin{COR*}{} Under the conditions of the theorem~\ref{THM:6-9}
the integral transform $J_a$ is an GFT iff the Fourier transform  $\widetilde a(\xi,t)$ with respect to $s$ of the kernel $a(s,t)$  satisfies the equation
\begin{equation}{}\label{6-13}
|\widetilde a(\xi,t)|=\Cr_K^{-1/2} |\omega (\xi,t)|,\quad\text{где}\quad
\omega (\xi,t)=\DOD{(\xi u(t))}{t}.
\end{equation}
\end{COR*}

Examples of GFT are the integral transforms $J_a$ with the kernels 
$a(s,t),$ for which the Fourier transforms with respect to $s$ are functions $\widetilde a(\xi,t)$ of the form
\begin{equation}{}\label{6-33}
\widetilde a(\xi,t)=\Cr_K ^{-1/2} \omega (\xi,t)\prod_{p=1}^l
         \xi_p^{\lambda _p}\overline\xi_p^{-\lambda _p},
\qquad \qquad \lambda _p\in\CC .
\end{equation}

Let us describe the kernels  $a(s,t)$ explicitly. Since $\omega (\xi,t)$ is a homogeneous polynom of $\xi_1, \dots ,\xi_l$ of the degree $k,$ then
the function  $\widetilde a(\xi,t)$ can be presented in the form
$$
\widetilde a(\xi,t)= \quad \sum_{m_1+ \ldots +m_l=k}
  \Bigl[ u_{m_1, \dots ,m_l} (t)
\prod_{p=1}^l (\xi_p^{m_p+ \lambda_p }\overline \xi_p
^{- \lambda_p })\Bigr].
$$
Therefore
$$
a(s,t)=\sum_{m_1+ \ldots +m_l=k}
\Bigl[ u_{m_1, \dots ,m_l} (t)
\prod_{p=1}^l \mathcal F(\xi_p^{m_p+ \lambda_p }\overline \xi_p^{- \lambda_p })\Bigr],
$$
where $\mathcal F$ is the inverse Fourier transform. In particular,
(see \cite{GGV}), if $\lambda _p\ne 0,$ $p=1, \dots ,l,$
then $\mathcal F(\xi_p^{m_p+ \lambda_p }\overline \xi_p^{- \lambda_p })$ is, up to a factor, the fuction 
$s_p^{- m_p- \lambda_p-1}\overline s_p^{\lambda _p-1}.$

In the special case $\lambda  _1= \ldots = \lambda _l=0$  the GFT is local and its kernel $a(s,t)$ has the form:
$$
a(s,t)=\Cr_K ^{-1/2} \omega (\DOD{}{s};t)\delta (s).
$$

{\bf 6. GFT related to the complex of $k$-dimensional planes in $\RR^n$}.
Consider a real analog of the complex of  $k$-dimensional planes.
Planes of the corresponding complex $K$ are defined by the same equations  \eqref{6-421}, where
$\alpha =(\alpha _1, \dots ,\alpha _l)\in\RR^l$ and
$t=(t_1, \dots ,t_k)\in\RR^k$, and generalized integral transforms 
$f \mapsto \phi=J_af$ are defined by the formula \eqref{6-52}, where $a(s,t)$ is a generalized function on $\RR^n$ and integration is performed respect to the Lebesgue measure on $\RR^n$.
The same relation  \eqref{6-59} holds between the Fourier transform of
$f$ and the Fourier transform of $\phi=J_af.$

Unlike the complex case, the Crofton function for the real complex $K$ is not constant in general  and it can be zero on an open set. The transform $J_a$ can be invertible iff this function is nonzero almost everywhere.

\begin{THM}\label{THM:66-9}
If the Crofton function  $\Cr_K(\eta,\xi)$ for the real $K$ is finite and nonzero almost everywhere and the kernel  $\widetilde a(\xi,t)$ is a usual function which is nonzero almost everywhere, then the following inversion formula holds for the integral transform $f \mapsto \phi=J_af:$
\begin{equation}{}\label{66-55}
f(x,y)=\int \phi(\alpha ,t ) A(y-u(t)x- \alpha ,t) d\mu(\alpha ) d\mu(t),
\end{equation}
where
$
\displaystyle
A(s,t)= \int \frac{|\omega (\xi,t)|}
{\Cr_K(-\xi u(t),\xi) \widetilde a(\xi,t)}\,
e^{i {\langle s,\xi \rangle}} d\mu(\xi).
$
\end{THM}
\begin{COR*}{} The integral transform $J_a$ is an GFT if the Fourier transform  $\widetilde a(\xi,t)$ with respect to $s$ of the kernel $a(s,t)$ satisfies the equation 
\begin{equation}{}\label{7-1}
|\widetilde a (\xi,t)|=\Cr_K^{-1/2} (\xi u(t),\xi)
        |\omega (\xi,t)|^{1/2}.
\end{equation}
\end{COR*}

In particular,  if all the elements of the matrix $u(t)$ are linear functions, then the Crofton function is identically equal to unity, hence $|\widetilde a (\xi,t)|=|\omega (\xi,t)|^{1/2}.$

\medskip

{\it Example. } In the case $k=1,$ $\Cr_K (\eta,\xi)\equiv 1$ and
$\widetilde a(\xi.t)=|\sum_{i=1}^{n-1} u'_i(t)\xi_i|^{1/2+i\rho},$
$\rho\in\RR$
the GFT kernel has the form:
$$
a(s;t)= \Bigl|\frac{s_1}{u'_1(t)}\Bigr|^{-3/2-i\rho}
\prod_{i=2}^{n-1} \delta  (u'_1(t) s_i-u'_i(t)s_1).
$$


\end{document}